\documentclass[twoside,12pt]{amsart} 
\usepackage[dvips]{graphicx}
\DeclareGraphicsExtensions{ps,eps}
\makeatletter
\def\serieslogo@{}
\makeatother
\makeatletter 
\def\@setcopyright{}
\makeatother

\newtheorem{Theorem}{Theorem}[section]
\newtheorem{Lemma}[Theorem]{Lemma}
\newtheorem{Definition}[Theorem]{Definition}

\theoremstyle{definition}

\theoremstyle{remark}

\newcommand{\ERRE}{{{\mathbb{R}}}}
\newcommand{\ZED}{{{\mathbb{Z}}}}
\newcommand{\rajtaro}{{{\rightarrow}}}
\newcommand{\LERAY}{{{\mathbb{P}}}}
\newcommand{\LDIV}{L_{\mathrm{div}}^2(D;{{\mathbb{R}}}^d)}
\newcommand{\LDIVo}{L^2_{{\mathrm{div}},0}(D;{{\mathbb{R}}}^d)}
\newcommand{\e}{\epsilon}

\begin{document}

\title
{
Topology-preserving diffusion of divergence-free vector fields
and magnetic relaxation
}

\author{Yann Brenier}
\address{
CNRS, Centre de math\'ematiques Laurent Schwartz,
Ecole Polytechnique, FR-91128 Palaiseau, France.
}
\curraddr{}
\email{brenier@math.polytechnique.fr}

\maketitle
\markboth{}{}

%\section*{PRELIMINARY VERSION 2 -PLEASE DO NOT CIRCULATE}
%%%%%%%%%%%%%%%%%%%%%%%%%%%%%%%%%%%%
\section*{Abstract}
%%%%%%%%%%%%%%%%%%%%%%%%%%%%%%%%%%%%

The usual heat equation is not suitable to preserve the topology of
divergence-free vector
fields, because it destroys their integral line structure. On the
contrary, in the
fluid mechanics literature, on can find
examples of topology-preserving diffusion
equations for divergence-free vector fields. They are very
degenerate since they admit all stationary solutions to the Euler
equations of incompressible fluids as equilibrium points.
For them, we provide a suitable concept of "dissipative solutions", 
which shares common features with both P.-L. Lions' dissipative
solutions to the Euler equations
and the concept of "curves of maximal slopes", \`a la De Giorgi,
recently used to study the scalar heat equation in very general
metric spaces.
We show that the initial value problem admits such global solutions,
at least in the two space variable case, 
and they are unique whenever they are smooth.

%%%%%%%%%%%%%%%%%%%%%%%%%%%%%%%%%%%%
\section{Introduction}
%%%%%%%%%%%%%%%%%%%%%%%%%%%%%%%%%%%%
Related to the numerous literature devoted to 
``topological fluid mechanics'' (see \cite{AK,FH,Mo,Ni,Ni2,Ri} and
many others),
there are some highly non-linear (and degenerate)
diffusion equations for divergence-free vector fields. 
A typical example is 
\begin{equation}
\label{MR}
\partial_t B+\nabla\cdot(B\otimes v-v\otimes B)=0,\;\;\;
v=\LERAY\nabla\cdot(B\otimes B),\;\;\;\nabla\cdot B=0,
\end{equation}
where $\LERAY$ denotes the $L^2$ projection onto divergence-free vector
fields.
Following \cite{Mo}, we call them ``magnetic relaxation equations''(MRE).
They have a (somewhat artificial) physical interpretation
as they (tentatively) describe a friction dominated model of incompressible  
magnetohydrodynamics (MHD) (in some lose sense, ``MHD in porous media''). 
%(A variant is
%$$
%-\Delta v=\LERAY\nabla\cdot(B\otimes B),
%$$
%and will be briefly discussed at the end of the paper.)
But, they have a specific mathematical interest because
of their link with the Euler equations of incompressible fluids and the theory
of "topological hydrodynamics", as discussed now.
\\
\\
Let us consider the MRE (\ref{MR}) 
on the flat torus $D=\ERRE^d/\ZED^d$, just for
simplicity, and sketch their three main properties.
First, these equations admits an interesting 
dissipation property for the ``magnetic energy'', namely
$$
\frac{d}{dt}\int \frac{|B|^2}{2} dx+\int |\LERAY\nabla\cdot(B\otimes B)|^2 dx=0
$$
(this is, of course, formally obtained by elementary calculations,
assuming $B$ to be smooth).
Next, the ``equilibrium states'', for which the energy no longer
dissipates, are precisely all stationary solutions to the 
Euler equations of homogeneous 
incompressible fluids, namely
$$
\LERAY\nabla\cdot(B\otimes B)=0,\;\;\;\nabla\cdot B=0.
$$
Last, equations (\ref{MR})
are "topology-preserving" in the sense that there is a velocity 
field $v=v(t,x)\in \ERRE^d$ that "transports" $B$,
i.e.
\begin{equation}
\label{TE}
\partial_t B+\nabla\cdot(B\otimes v-v\otimes B)=0.
\end{equation}
What we mean by "topology-preserving"
is that, for each fixed time $t$, the "magnetic lines" of
$B$ (i.e. the integral lines $s\rightarrow \xi(s)$ satisfying
$\xi'(s)=B(t,\xi(s))$, where $t$ is ``frozen'')
are ``transported'' by the flow associated to $v$, as time evolves.
Thus, these lines
keep their topology unchanged
during the evolution, in particular their knot structure.
[More precisely, let us use "material coordinates" 
$(t,a)\rightarrow X(t,a)$, so that
$$
\partial_t X(t,a)=v(t,X(t,a)),\;X(0,a)=a.
$$
Then, under suitable smoothness assumptions on $v$ and $B$, 
the transport equation (\ref{TE}) exactly means
$
B_i(t,X(t,a))=\sum_j \partial_jX_i(t,a)B_j(0,a),\;\;i=1,\cdot\cdot\cdot,d
$
(to check the formula, just differentiate it in $t$ and use the chain rule.)
This implies that the magnetic lines of $B$ at time $t$
are the images by $X(t,\cdot)$ of those of $B$ at times 0.
(To check this statement, differentiate $X(t,\xi(s))$ in $s$
for each integral line $\xi$ of $B$ at time 0.)]
\\
To summarize these three properties, we can say that the ``magnetic relaxation
equations" (\ref{MR}) are a good candidate to drive,
as time goes to infinity, each given initial magnetic
field of prescribed topology 
to a stationary solution of the Euler equations with the same topology.
This is clearly part of the more ambitious 
program of "topological hydrodynamics", 
as developed
in the papers of Moffatt \cite{Mo} and the book of Arnold and Khesin \cite{AK}.
\\
\\
Let us emphasize that the standard linear
diffusion equation (on the flat-torus) for  divergence-free vector fields reads
$
\;\;\partial_t B=\Delta B\;\;
$
and is certainly not ``topology-preserving''
since it cannot be 
written as a transport equation (\ref{TE})
for any vector field $v=v(t,x)$.
This is in sharp contrast with the standard linear heat equation for positive 
density fields
$
\partial_t \rho=\Delta\rho,
$
which can be easily put in ``transport'' form
\begin{equation}
\label{heat}
\partial_t \rho+\nabla\cdot(\rho v)=0,\;\;\;v=-\nabla(\log\rho).
\end{equation}
Recently, the heat equation has been studied by Gigli, Gigli-Kuwada-Ohta,
Ambrosio-Gigli-Savar\'e
\cite{G,GKO,AGS2}, in a very general class of metric spaces. Their method
is based on the following 
very simple and remarkable idea (that combines the concept of "curves of
maximal slopes" introduced
by the De Giorgi school for "gradient flows'' and the interpretation by Kinderlehrer,
Jordan and Otto \cite{JKO,AGS} of the heat equation as the gradient flow
of Boltzmann's entropy for a suitable Monge-Kantorovich metric
on the set of probability measures).
We first say that a pair of measures $(\rho(t,x)\ge 0,q(t,x)\in \ERRE^d)$ 
is admissible if it
solves the ``continuity equation''
\begin{equation}
\label{CE}
\partial_t \rho+\nabla\cdot q=0.
\end{equation}
Next, we formally get, for each admissible pair
$$
\frac{d}{dt}\int 2\rho\log\rho=2\int 
\frac{\nabla \rho\cdot q}{\rho}
$$
$$
=\int \frac{|q+\nabla\rho|^2}{\rho}
-\int \frac{|q|^2}{\rho}
-\int \frac{|\nabla\rho|^2}{\rho}\;.
$$
This suggests to characterize the solutions of the heat equation precisely
as those admissible pairs $(\rho,q)$ that achieve inequality
\begin{equation}
\label{HV}
2\frac{d}{dt}\int \rho\log\rho 
+\int \frac{|q|^2}{\rho}
+\int \frac{|\nabla\rho|^2}{\rho}\le 0\;.
\end{equation}
This very simple formulation is quite powerful. First, the set of solutions
$(\rho,q)$ for a given initial condition is convex and weakly compact.
(Notice that the three functionals involved in (\ref{HV}) are all convex 
in $(\rho,q)$,
with possible value $+\infty$, the first one being strictly convex.)
Next, uniqueness of solutions directly
follows from the strict convexity of $\rho\rightarrow \rho\log\rho$.
(Indeed, the average of two distinct solutions would lead to a strict inequality
in (\ref{HV}), which turns out to be not possible). Finally, formulation (\ref{HV})
makes sense in a very large class of metric spaces \cite{G,GKO,AGS2}. 
Notice
that this strategy can also be carried out for a rather large class of 
non-linear diffusion equations, as explained in Appendix 1.
\\
\\
We have spent several lines in explaining this recent approach
to the scalar heat equation precisely because we are going to follow a 
similar (but much
less successful) way to address
the (more) challenging magnetic relaxation equations
(\ref{MR}).
First, we call admissible solution any pair of time-dependent 
divergence-free vector fields $(B,v)$ solving the transport equation 
(\ref{TE}), namely
\begin{equation}
\label{TEter}
\partial_t B+\nabla\cdot(B\otimes v-v\otimes B)=0
\end{equation}
(which, unfortunately, is non-linear but, at least, has a nice ``div-curl''
structure in the two-dimensional case $d=2$).
Next, we observe that, for any (smooth) admissible pair $(B,v)$:
$$
\frac{d}{dt}\int |B|^2 dx=2\int v\cdot \nabla\cdot(B\otimes B)dx.
$$
(just because of (\ref{TE}))
$$
=2\int v\cdot \LERAY\nabla\cdot(B\otimes B)dx
$$
(because $v$ is divergence-free)
$$
=\int|v-\LERAY\nabla\cdot(B\otimes B)|^2dx
-\int|v|^2dx
-\int|\LERAY\nabla\cdot(B\otimes B)|^2dx.
$$
Therefore, we
can characterize the solutions $(B,v)$ of (\ref{MR}),
just by asking them to be admissible and satisfy the following inequality
$$
\frac{d}{dt}\int |B|^2 dx+\int|v|^2dx
+\int|\LERAY\nabla\cdot(B\otimes B)|^2dx\le 0.
$$
Unfortunately, with respect to the simpler case of the scalar heat equation,
we loose on two sides. First, the transport ``constraint'' (\ref{TEter})
is not linear (in contrast with the continuity equation (\ref{CE})).
Second, the energy inequality involves a non-convex functional of $B$,
namely $\int|\LERAY\nabla\cdot(B\otimes B)|^2dx$.
However, in the present paper, we will be able to overcome some of
these difficulties. The output is a ``reasonable''
concept of ``dissipative solutions'', sharing the same strength and weakness as Lions'
dissipative solutions to the Euler equations \cite{Li}: weak compactness
(at least in the two-dimensional case $d=2$, in our case)
and uniqueness whenever they are smooth. To finish this introduction,
let us mention the analysis by Nishiyama of related magnetic relaxation equations 
\cite{Ni,Ni2}, based on the concept
of measure-valued solutions, as well as two recent papers 
\cite{BDS,DST} on the concept of dissipative
measure-valued solutions for elastodynamics and fluid mechanics.

%%%%%%%%%%%%%%%%%%%%%%%%%%%%%%%%%%%%
\section{Dissipative solutions to the magnetic relaxation equations}
%%%%%%%%%%%%%%%%%%%%%%%%%%%%%%%%%%%%

\subsection*{Preliminaries}
%%%%%%%%%%%%%%%%%%%%%%%%%%%%%%%%%%%%

Let $D$ be the $d$-dimensional
flat torus $(\ERRE/\ZED)^d$, and $||\cdot||$ be the $L^2$-norm on $D$.
Let us denote $\LDIVo$ the space of all
$L^2$ vector fields on $D$ that are divergence-free 
%and parallel to the boundary.
%(Observe that  $\LDIVo$ may be more precisely defined as the $L^2$ closure of all smooth 
%compactly supported divergence-free vector fields
%on $D$.) 
\\
Next we define 
\begin{equation}
\label{Leray}
L(B)=||\LERAY\nabla\cdot(B\otimes B)||^2
\in [0,+\infty],
\end{equation}
for all $B\in \LDIVo$,
where $\LERAY$ denotes the (Helmholtz-Leray) orthogonal projector 
$L^2\rajtaro \LDIVo$.
A more precise definition is obtained by duality
\begin{equation}
\label{Leray2}
\begin{split}
L(B)=\sup\int_D (B\otimes B):(\nabla z+\nabla z^T)dx -||z||^2 \;,
\\
z\in C^1(D\;;\;\ERRE^d),
\;\;\;\nabla\cdot z=0,
\end{split}
\end{equation}
where 
$(B\otimes B):(\nabla z+\nabla z^T)$ should be understood as
$$
\sum_{i,j=1}^d B_i B_j(\partial_i z_j+\partial_j z_i).
$$
 $L$ is not a convex function of $B$ and we have to find a substitute for
$L$. This is why,
for each real nonnegative number $r\in \ERRE_+$, we define, for  $B\in \LDIV$
\begin{equation}
\label{KRB}
\begin{split}
K_r(B)=\sup\{\int_D (B\otimes B):(\nabla z+\nabla z^T+rI)dx-||z||^2\;,
\\
z\in C^1(D\;;\;\ERRE^d),
\;\;\;\nabla\cdot z=0,
\;\;\;\nabla z+\nabla z^T+rI\ge 0\}
\end{split}
\end{equation}
(in the sense of symmetric matrices, where $I$ denotes the identity matrix).
Notice that $K_r(B)$ is always bounded from below by $r||B||^2$ (take $z=0$ in its definition)
and is a convex function of $B$ (as a supremum of
positively curved quadratic functions of $B$, thanks to the constraint
$\nabla z+\nabla z^T+rI\ge 0$).
In addition, we can recover $L$ out of the $K_r$ since 
%\begin{equation}
%\begin{split}
%L(B)=
%\sup\{\int_D (B\otimes B):(\nabla z+\nabla z^T)dx-||z||^2,
%\\
%\;\;z\in C^1(D\;;\;\ERRE^d),
%\;\;\;\nabla\cdot z=0,
%\;\;\;\nabla z+\nabla z^T+rI\ge 0,\;\;\;r\ge 0.\}
%\\
%=\sup_{r\ge 0}\;K_r(B)-r||B||^2.
%\end{split}
%\end{equation}

\begin{equation}
\begin{split}
\sup_{r\ge 0}\;K_r(B)-r||B||^2
=\sup\{\int_D (B\otimes B):(\nabla z+\nabla z^T+rI)dx-||z||^2-r||B||^2,
\\
\;\;r\ge 0,\;\;\;z\in C^1(D\;;\;\ERRE^d),
\;\;\;\nabla\cdot z=0,
\;\;\;\nabla z+\nabla z^T+rI\ge 0\}
\\
=\sup\{\int_D (B\otimes B):(\nabla z+\nabla z^T)dx-||z||^2
\\
\;\;r\ge 0,\;\;\;z\in C^1(D\;;\;\ERRE^d),
\;\;\;\nabla\cdot z=0,
\;\;\;\nabla z+\nabla z^T+rI\ge 0\}
\\
=L(B).
\end{split}
\end{equation}

\subsection*{Definition of dissipative solutions to the magnetic relaxation equations}
%%%%%%%%%%%%%%%%%%%%%%%%%%%%%%%%%%%%

\begin{Definition}
Given a final time $T>0$ and $B^0\in \LDIVo$, we say that a pair 
$$(t,x)\in [0,T]\times D\rajtaro (B,v)(t,x)\in \ERRE^d\times \ERRE^d$$
is a dissipative solution of the
MRE (\ref{MR}) with initial condition $B^0$, if
\\
i) $B$ is weakly continuous from $[0,T]$ to $\LDIVo$ with $B(0)=B^0$;
\\
ii) $v$ is square-integrable  from $[0,T]$ to $\LDIVo$;
\\
iii) $B,v$ solves the transport equation (\ref{TE}), namely
$$
\partial_t B+\nabla\cdot(B\otimes v-v\otimes B)=0,
$$
in the sense of distributions;
%\begin{equation}
%\label{TEbis}
%\begin{split}
%\frac{d}{dt}\int_D B\cdot\beta
%=\int_D B_{ti}(\beta_{ti,t}+v_{tj}(\beta_{ti,j}-\beta_{tj,i}))
%\\
%{\mathrm{a.e.}} \;t\in [0,T], \;\;\forall \beta\in C^1(D\;;\;\ERRE^d)
%\end{split}
%\end{equation}
%(where we use notations $\beta_{ti,j}=\partial_j\beta_{ti}$, etc...and skip
%summations on repeated indices $i$, $j$,...);
\\
iv) 
\begin{equation}
\label{entropy}
\begin{split}
||B(t,\cdot)||^2+\int_0^t [||v(s,\cdot)||^2+K_{r(s)}(B(s,\cdot))]\exp(R(t)-R(s))ds
\\
\le||B(0,\cdot)||^2\exp(R(t))
\;\;\;\forall t\in[0,T],
\end{split}
\end{equation}
for all nonnegative function $t\rightarrow r(t)\ge 0$, with $R(t)=\int_0^t r(s)ds$,
where $K_r$ is defined by (\ref{KRB}).
\end{Definition}

%%%%%%%%%%%%%%%%%%%%%%%%%%%%%%%%%%%%
\section{Stability of smooth solutions among dissipative solutions}
%%%%%%%%%%%%%%%%%%%%%%%%%%%%%%%%%%%%
We ignore whether or not the MRE (\ref{MR}) are locally well-posed
in any space of smooth functions. (As a matter of fact, this is a very
interesting issue, since these equations can be considered
as parabolic only in a weak sense.)
However, we can prove:

\begin{Theorem}
\label{uniqueness}
Assume $D=(\ERRE/\ZED)^d$.
Let $(B,v)$ and $(\beta,\omega)$ be respectively a dissipative and a smooth solution to
the MRE (\ref{MR}) up to time $T$. Then, there is a constant $C$ depending only on the
spatial Lipschitz constant of $(\beta,\omega)$, so that, for all $t\in [0,T]$,
\begin{equation}
\label{strong uniqueness}
\begin{split}
||(B-\beta)(t,\cdot)||^2
+ \int_0^t \exp(C(t-s))\frac{1}{2}||v(s,\cdot)-\omega(s,\cdot)||^2ds\;\le 
||(B-\beta)(0,\cdot)||^2\exp(Ct)
\end{split}
\end{equation}
\end{Theorem}
This implies the uniqueness of smooth
solutions among all dissipative solutions, for any given prescribed 
smooth initial condition.
\subsection*{Proof}
For simplicity, we use notations $B_t$, $v_t$ for $B(t,\cdot)$, $v(t,\cdot)$, etc...
Let us introduce for each $t\in [0,T]$
$$
N_t=||B_0||^2\exp(rt)-\int_0^t (||v_s||^2+K_r(B_s))\exp(r(t-s))ds
$$
where $r$ is a nonnegative constant to be chosen later. By definition, we have
$$
(\frac{d}{dt}-r)N_t+||v_t||^2+K_r(B_t)=0
$$
(in the distributional sense and also for a.e. $t\in [0,T]$).
Since $(B,v)$ is a dissipative solution we get from (\ref{TE}) 
$$
N_t\ge ||B_t||^2, \forall t\in [0,T].
$$
We now set
$$
e(t)=
||B_t-\beta_t||^2+(N_t-||B_t||^2)
=N_t-2((B_t,\beta_t))+||\beta_t||^2,\;\forall t\in [0,T].
$$
where $((\cdot,\cdot))$ denotes the $L^2$ inner product,
and compute the time derivative $e'(t)$ of $e(t)$.
We already know that
$$
\frac{d}{dt}N_t=rN_t-||v_t||^2-K_r(B_t)
$$
Next, since $(B,v)$ is a dissipative solution we get from (\ref{TE})
$$
\frac{d}{dt}((B_t,\beta_t))
=\int B_{ti}(\beta_{ti,t}+v_{tj}(\beta_{ti,j}-\beta_{tj,i}))
$$
where we use notations $\beta_{ti,j}=\partial_j(\beta_t)_{i}$, etc...and skip
summations on repeated indices $i$, $j$...
Thus, we find
$$
e'(t)=rN_t-||v_t||^2-K_r(B_t)+
\int 2(\beta-B)_{ti}\beta_{ti,t}-2B_{ti}v_{tj}(\beta_{ti,j}-\beta_{tj,i}).
$$
By definition (\ref{KRB}), we can find a constant
$r$, depending only on the spatial Lipschitz constant of $\omega$,
large enough so that 
(by setting $z=-\omega$ in (\ref{KRB}))
$$
-K_r(B_t)\le \int B_{ti}B_{tj}(\omega_{ti,j}+\omega_{tj,i})dx
-r||B_t||^2+||\omega_t||^2.
$$
Thus
$$
e'(t)\le r(N_t-||B_t||^2)-||v_t||^2+||\omega_t||^2+J_t
$$
where
$$
J_t=\int B_{ti}B_{tj}(\omega_{ti,j}+\omega_{tj,i})
+2(\beta-B)_{ti}\beta_{ti,t}-2B_{ti}v_{tj}(\beta_{ti,j}-\beta_{tj,i}).
$$
We may write
$$
J_t=
J_t^{Q}+J_t^{L1}+J_t^{L2}+J_t^{C}
$$
where 
$J_t^{Q}$, $J_t^{L1}$, $J_t^{L2}$, $J_t^{C}$ are respectively 
quadratic, linear, linear, and constant with respect to $B-\beta$ and $v-\omega$,
with coefficient depending only on $\omega,\beta$:
$$
J_t^{Q}=\int (B-\beta)_{ti}(B-\beta)_{tj}(\omega_{ti,j}-\omega_{tj,i})
-2(B-\beta)_{ti}(v-\omega)_{tj}(\beta_{ti,j}-\beta_{tj,i})
$$
$$
J_t^{L1}=\int 2(B-\beta)_{ti}[\beta_{tj}(\omega_{ti,j}+\omega_{tj,i})-\beta_{ti,t}
-\omega_{tj}(\beta_{ti,j}-\beta_{tj,i})]
$$
$$
J_t^{L2}=\int 2(v-\omega)_{tj}\beta_{ti}\beta_{tj,i}
$$
$$
J_t^{C}=\int [\beta_{ti}\beta_{tj}(\omega_{ti,j}+\omega_{tj,i})
-2\beta_{ti}\omega_{tj}(\beta_{ti,j}-\beta_{tj,i})].
$$
Let us reorganize these four terms.
By integration by part of its second term, we see that $J_t^{C}=0$, using
that $\beta$ and $\omega$ are divergence-free.
Next, since $B_t-\beta_t$ is divergence-free, we have
$$
J_t^{L1}=\int 2(B-\beta)_{ti}[
-\beta_{ti,t}-\omega_{tj}\beta_{ti,j}+\beta_{tj}\omega_{ti,j}]
$$
and we may reorganize 
$$
J_t^{L2}=2((v_t-\omega_t,\omega_t))+
\int 2(v-\omega)_{tj}[\beta_{ti}\beta_{tj,i}-\omega_{tj}]
$$
Thus
$$
e'(t)\le r(N_t-||B_t||^2)-||v_t-\omega_t||^2+J_t^{Q}+J_t^{L}
$$
where
$$
J_t^{L}=-2((B_t-\beta_t,\partial_t \beta_{t}+(\omega_{t}\cdot\nabla)\beta_{t}
-(\beta_{t}\cdot\nabla)\omega_{t}))
+2((v_t-\omega_t,\nabla(\beta_{t}\otimes\beta_{t})-\omega_{t})).
$$
Clearly
$$
J_t^{Q}\le \frac{1}{2}||v_t-\omega_t||^2+C'||B_t-\beta)_t||^2
$$
for some constant $C'$ depending only on the spatial 
Lipschitz constant of $(\beta,\omega)$. Finally,
we have obtained
$$
e'(t)+\frac{1}{2}||v_t-\omega_t||^2\le r(N_t-||B_t||^2)+C'||B_t-\beta_t||^2+J_t^{L},
$$
and, therefore, by definition of $e(t)$
$$
e(t)=
||B_t-\beta_t||^2+(N_t-||B_t||^2),
$$
we get
$$
e'(t)+\frac{1}{2}||v_t-\omega_t||^2\le Ce(t)+J_t^{L},
$$
where $C$ is another constant depending only on the spatial Lipschitz
constants of $(\beta,\omega)$ (through $R$).
By integration we deduce
$$
e(t)+\int_0^t e^{(t-s)C}\frac{1}{2}||v_s-\omega_s||^2ds
\le e(0)e^{tC}+\int_0^t e^{(t-s)C}J_s^{L}ds.
$$
Finally, let us remember that
$e(t)\ge ||B_t-\beta_t||^2$ with equality for $t=0$
(since $N_t\ge ||B_t||^2$ with equality at $t=0$).
Thus, we have shown
\begin{Lemma}
Assume $D=(\ERRE/\ZED)^d$. Let us fix $T>0$.
Let $(B,v)$ be a dissipative solution of the MRE (\ref{MR}) up to time $T$,
and $(\beta,\omega)$ be any pair of smooth functions $(\beta,\omega)$ chosen
in $\LDIVo$. Then there is a constant $C$ depending only on 
the spatial Lipschitz constant of $(\beta,\omega)$, up to time $T$,
so that, for all $t\in [0,T]$,
\begin{equation}
\label{strong uniqueness2}
\begin{split}
||B_t-\beta_t||^2+\int_0^t e^{(t-s)C}\frac{1}{2}||v_s-\omega_s||^2ds
\le ||B_0-\beta_0||^2e^{tC}+\int_0^t e^{(t-s)C}J_s^{L}ds,
\\
J_t^{L}=-2((B_t-\beta_t,\partial_t \beta_{t}+(\omega_{t}\cdot\nabla)\beta_{t}
-(\beta_{t}\cdot\nabla)\omega_{t}))
+2((v_t-\omega_t,\nabla(\beta_{t}\otimes\beta_{t})-\omega_{t})).
\end{split}
\end{equation}
\end{Lemma}
We now see that $J_t^L$ exactly vanishes as
$(\beta,\omega)$ is a smooth solution to the MRE (\ref{MR}).
$$
\partial_t\beta+(\omega\cdot\nabla)\beta-(\beta\cdot\nabla)\omega=0,
\;\;\;\omega=\LERAY\nabla(\beta\otimes\beta),\;\;\;
\nabla\cdot\beta=\nabla\cdot\omega=0,
$$
and the proof of Theorem \ref{uniqueness} immediately follows.
%%%%%%%%%%%%%%%%%%%%%%%%%%%%%%%%%%%%
\section{Existence of dissipative solutions in two space dimensions}
%%%%%%%%%%%%%%%%%%%%%%%%%%%%%%%%%%%%

\begin{Theorem}
\label{existence}
Assume $d=2$ and $D=(\ERRE/\ZED)^2$. 
Let $T>0$ and fix an initial condition $B^0\in \LDIVo$.
Then there is at least one dissipative solution $(B,v)$ of
the MRE (\ref{MR}) up to time $T$. 
This solution can be obtained as the limit in
$C^0([0,T];\LDIVo_w)\times L^2([0,T];\LDIVo)_w$,
as parameters $(\e,\mu,\nu)$ go to zero,
of the unique solution of
the MHD system (with friction and viscosity)
\begin{equation}
\label{MHD}
\begin{split}
\e(\partial_t v+\nabla\cdot v\otimes v)+v+\nabla p
-\mu\Delta v=\nabla\cdot B\otimes B,
\\
\partial_t B+\nabla\cdot(B\otimes v-v\otimes B)-\nu\Delta B=0,
\;\;\;\;\nabla\cdot v=\nabla\cdot B=0,
\end{split}
\end{equation}
with smooth initial conditions chosen so that
$B(0,\cdot)$, $\sqrt{\e}\;v(0,\cdot)$ approach respectively $B^0$ and $0$ in $L^2$.
\end{Theorem}

\subsection*{Proof}
Since $d=2$, the MHD system has global smooth solutions for any
smooth intial condition \cite{ST}.
From these MHD equations, respectively multiplied by $v$ and $B$ and integrated over
$D$, we get two straightforward estimates
\begin{equation}
\label{MHD estimates}
\begin{split}
\e\frac{d}{dt}\int \frac{|v|^2}{2}+\int |v|^2+\mu\int |\nabla v|^2
=-\int (B\otimes B):\nabla v
\\
\frac{d}{dt}\int \frac{|B|^2}{2}+\nu\int |\nabla B|^2
=\int (B\otimes B):\nabla v
\end{split}
\end{equation}
We first add up these estimates in order to get the total energy
balance
\begin{equation}
\label{MHD balance}
\frac{d}{dt}\int \frac{|B|^2+\e|v|^2}{2}
+\int |v|^2+\mu\int |\nabla v|^2+\nu\int |\nabla B|^2
=0.
\end{equation}
This implies that $B$ and $v$ are respectively uniformly bounded
(with respect to $(\e,\mu,\nu)$)
in $L^\infty([0,T],L^2(D))$ and $L^2([0,T],L^2(D))$.
Using  the second equation of (\ref{MHD}),
we also see that $B$ is uniformly bounded in
$C^{1/2}([0,T],({\rm{Lip}}(D;\ERRE^d)')$ (where ${\rm{Lip}}(D;\ERRE^d)'$ denotes
the dual of the space of vector-valued Lipschitz functions). Indeed:
\begin{equation}
\begin{split}
\label{aubin}
\forall t,s\in [0,T],\;\;\;|\int (B(t,x)-B(s,x))\cdot z(x)dx|^2=
\\
|\int_s^t dt'\int [(B(t',x)\otimes v(t',x)):(\nabla z^T-\nabla z)(x)
-\nu \nabla B(t',x):\nabla z(x)]dx|^2
\\
\le |t-s|\;{\rm{Lip(z)}}^2
\;\;[4\sup_{t"\in [0,T]}\int |B(t",x)|^2 dx \int |v(t',x)|^2dxdt'
+\nu^2\int |\nabla B(t',x)|^2dxdt'].
\end{split}
\end{equation}
We deduce that, as $\e,\mu,\nu$ go to zero,
$B$ and $v$ are compact respectively in the spaces
$C^0([0,T],\LDIVo_w)$ 
%(thanks to Aubin's lemma, cf.  \cite{Li})
and $L^2([0,T],\LDIVo)_w$, where subscript $w$
refers to the weak topology of $L^2$.
\\
Next, we use the weak formulation (in $x$) of the first equation of (\ref{MHD}),
namely:
$$
\e(\partial_t v+\nabla\cdot v\otimes v)+v+\nabla p
-\mu\Delta v=\nabla\cdot B\otimes B,
$$
to get, for any fixed smooth test function $z(t,\cdot)$ valued in $\LDIVo$
$$
\e\frac{d}{dt}\int  v_i z_i
-\int \e (v_i z_{i,t}+v_i v_{j} z_{i,j})
+\int (z_i v_i 
+\mu z_{i,j} v_{i,j}+B_{i} B_j z_{i,j})=0
$$
Adding up the energy balance (\ref{MHD balance}), namely:
$$
\frac{d}{dt}\int (\frac{|B|^2+\e|v|^2}{2})
+\int |v|^2+\mu\int |\nabla v|^2+\nu\int |\nabla B|^2
=0,
$$
we obtain
\begin{equation}
\label{MHD dissipation}
\begin{split}
\frac{d}{dt}\int (\frac{|B|^2
+\e|v|^2}{2}+\e v_i z_i)
-\int (\e(v_i z_{i,t}+v_i v_{j} z_{i,j})+\mu z_{i,j} v_{i,j})
\\
+\int (\frac{|v|^2}{2}+B_{i} B_jz_{i,j}-\frac{|z|^2}{2})
+\int(\frac{|v+z|^2}{2}+\mu|\nabla v|^2+\nu |\nabla B|^2)=0.
\end{split}
\end{equation}
Let us now introduce any nonnegative function $t\rightarrow r(t)$ and $R(t)=\int_0^t r(s)ds$,
and assume that $(z_{i,j}+z_{j,i}+r(t)\delta_{ij})$ is a positive matrix.
After multiplication by 2, we may rearrange
(\ref{MHD dissipation}) as
\begin{equation}
%\label{MHD dissipation2}
\begin{split}
(\frac{d}{dt}-r)\int (|B|^2+\e|v+z|^2)
+\int (|v|^2+B_{i} B_j(z_{i,j}+z_{j,i}+r\delta_{ij})-|z|^2)
\\
+r\int \e|v+z|^2
-\frac{d}{dt}\int \e|z|^2
-2\int (\e(v_i z_{i,t}+v_i v_{j} z_{i,j})+
+\mu z_{i,j} v_{i,j})
\\
+\int (|v+z|^2+2\mu|\nabla v|^2+2\nu|\nabla B|^2)=0.
\end{split}
\end{equation}
Thus,
\begin{equation}
\begin{split}
(\frac{d}{dt}-r)\int (|B|^2+\e|v+z|^2)
+\int (|v|^2+B_{i} B_j(z_{i,j}+z_{j,i}+r\delta_{ij})-|z|^2)
\le \eta(t)
\end{split}
\end{equation}
where $\eta(t)$ depends on the fixed test function $z$
and goes to zero in $L^1([0,T])$ with $(\e,\mu,\nu)$,
since $v$ is uniformly bounded in $L^2$.
Next, we integrate in time this differential inequality, and, then,
we let $(\e,\mu,\nu)$ go to zero, assuming that the initial condition
$B(0,\cdot)$ and $\sqrt{\e}v(0,\cdot)$ converge in $L^2$
respectively to the given initial condition $B^0$ and to $0$.
After these operations, we obtain for any accumulation point
of the $(B,v)$, still denoted by $(B,v)$,
\begin{equation}
%\label{entropy}
\begin{split}
||B(t,\cdot)||^2+\int_0^t ||v(s,\cdot)||^2\exp(R(t)-R(s))ds
\\
+\int_0^t \int_D dx(B_{i} B_j(2z_{i,j}+r\delta_{ij})-|z|^2)(s,x)\exp(R(t)-R(s))ds
\\
\le ||B(0,\cdot)||^2\exp(R(t))
\;\;\;\forall t\in[0,T],
\end{split}
\end{equation}
using the positivity of $(z_{i,j}+z_{j,i}+r\delta_{ij})$.
Next, taking the supremum with respect to $z$, for fixed $r$, leads
to the dissipation inequality (\ref{entropy}) 
involved in the dissipative formulation of the magnetic relaxation equations. 
However, we are still left with the problem of passing to the
limit in the transport equation (\ref{TE}). We do not see any clue for that,
except in the bidimensional case $d=2$, where we have a nice "div-curl" structure. 
Indeed, as $d=2$, at least locally, we can write
$$
B=(\partial_2 A,-\partial_1 A),
$$
for some scalar potential $A(t,x_1,x_2)\in\ERRE$.
Then, the transport equation (\ref{TE}) can be integrated out as
$$
\partial_t A+\nabla\cdot(Av)=0,
$$
and we can pass to the limit, since $|\nabla A|=|B|$ and $v$ are
well controled in $L^2$ (using estimate (\ref{aubin}) to handle the time
dependence). This concludes the proof.
%%%%%%%%%%%%%%%%%%%%%%%%%%%%%%%%%%%%%%%%%%%
%%%%%%%%%%%%%%
\section{
Appendix
}
\subsection{A general framework for dissipative equations}
A rather general framework that one can encounter in several situations of 
Mechanics and Physics
is as follows. (We do not claim any novelty in it, see, for instance, \cite{Mi}
for somewhat related issues.)
Working on the $d$-dimensional flat torus $D=\ERRE^d/\ZED^d$, for simplicity,
we call admissible a pair $(B,E)$ made of two time-dependent differential forms, one of degree $k$, say $B$,
that we assume to be closed $dB=0$,
and one of degree $k-1$, say $E$, linked by
\begin{equation}
\label{admissible}
\partial_t B+dE=0.
\end{equation}
(A simple example being $k=d$, as in the scalar heat equation where $B=\rho$, $E=q$.)
Next, we are given a scalar function $L(E,B)$, called "Lagrangian", that we suppose convex in $E$
(with possible value $+\infty$). Then, we define its Legendre-Fenchel
transform with respect to $E$, called "Hamiltonian'
\begin{equation}
\label{LF}
H(D,B)=\sup_E E\cdot D-L(E,B),
\end{equation}
and introduce the "defect" function
\begin{equation}
\begin{split}
\label{defect}
{\rm{def}}(B,E,D)=L(E,B)+H(D,B)-E\cdot D\ge 0
\\
{\rm{with\;equality\;if\;and\;only\;if\;}}E=\partial_1 H(D,B).
\\
{\rm{if\;and\;only\;if\;}}D=\partial_1 L(E,B).
\end{split}
\end{equation}
Now, we are given a convex function $\theta$
and compute for any (smooth) admissible pair $(E,B)$
$$
\frac{d}{dt}\int \theta(B)dx=-\int \theta'(B)\cdot dE
=-\int E\cdot D
$$
where $D=\delta(\theta'(B))$
(here $\delta=(-1)^{k-1} *^{-1}d *$ is the Hodge co-differential).
Thus
$$
\frac{d}{dt}\int \theta(B)dx
=\int ({\rm{def}}(B,E,D)-L(E,B)-H(D,B))\ge -\int (L(E,B)+H(D,B))
$$
with equality if and only if $E=\partial_1 H(D,B)$ (pointwise).
This suggests a "dissipative formulation" for the (highly) non-linear equation
\begin{equation}
\begin{split}
\label{general heat}
\partial_t B+dE=0,
\\
E=\partial_1 H(\delta(\theta'(B)),B)
\end{split}
\end{equation}
We call dissipative solutions of this equation any admissible pair $(E,B)$
such that
\begin{equation}
\label{general dissip}
\begin{split}
\frac{d}{dt}\int \theta(B)dx
+\int (L(E,B)+H(D,B))
\le 0,
\\
{\rm{where\;\;}}D=\delta(\theta'(B)).
\end{split}
\end{equation}
Of course, if 
$$
\int (L(E,B)+H(\delta(\theta'(B)),B))
$$
turns out to be a convex function of the pair $(E,B)$ the analysis gets rather simple,
but there is little chance to get interesting examples of this type unless $k=d$.
As a matter of fact, in the case $k=d$, denoting 
$(E,B)=(q,\rho)$, we get
the rather general non-linear scalar diffusion equation
$$
\partial_t\rho=\nabla\cdot \partial_1 H(\nabla(\theta'(\rho)),\rho).
$$
In the special case $L(q,\rho)=\rho c(\frac{q}{\rho})$, we get $H(v,\rho)=\rho c^*(v)$,
where $c$ is a convex function and $c^*$ its Legendre-Fenchel transform.
The resulting equation reads
$$
\partial_t\rho=\nabla\cdot (\rho \nabla c^*(\nabla(\theta'(\rho))).
$$
(Notice that this equation can be also handled by optimal transport methods \cite{Vi}, using "cost function" $(x,y)\rightarrow c(x-y)$ and "entropy function" $\theta$.)
The further choice 
$$
\theta(\rho)=\rho\log\rho,\;\;\;c(w)=|w|^2/2
$$
leads to the linear heat equation.
Another example
is the "relativistic heat equation" \cite{ACM,MP}, for which
$$
\theta(\rho)=\rho\log\rho,\;\;\;c(w)=-\sqrt{1-|w|^2}.
$$
\subsection{A ``topology-preserving'' diffusion equation for divergence-free vector fields, based on Born-Infeld Electromagnetism}
In order to provide a non-scalar application of 
the general framework, let us consider the special case $d=3$ and $k=2$. So $B$ is a closed $2-$ form in three space dimensions, which
corresponds to a divergence-free vector field, while $E$ is a $1-$ form, i.e. a vector field.
We use classical notations $\times$ for the wedge product
as well as $\nabla\times$ for $d$, which here is the curl operator.
So, we call admissible solutions any pair of fields $(E,B)(t,x)\in \ERRE^3$
satisfying
\begin{equation}
\label{BI1}
\partial_t B+\nabla\times E=0.
\end{equation}
Let us introduce the Born-Infeld "Lagrangian" \cite{BI}, parameterized by constant
$\lambda>0$:
\begin{equation}
\label{BI lag}
L_\lambda(E,B)=-\sqrt{\lambda^2+|B|^2-|E|^2-\lambda^{-2}(E\cdot B)^2}.
\end{equation}
Function $L_\lambda$, for each fixed value of $B$ is convex in $E$ (with infinite
value as the term under the square root gets negative). The Legendre-Fenchel
transform with respect to $E$ can be easily computed and is given by the
"Hamiltonian"
\begin{equation}
\begin{split}
\label{BI ham}
H_\lambda(D,B)=\sqrt{\lambda^2+|B|^2+\lambda^2|D|^2+|D\times B|^2}
\\
=\sqrt{(\lambda^2+|B|^2)(1+|D|^2)-(D\cdot B)^2}.
\end{split}
\end{equation}
Let us now consider a convex function $\theta:\ERRE^3\rightarrow \ERRE$.
From the general framework, we get a
"dissipative formulation" for the (highly) non-linear equation
\begin{equation}
\begin{split}
\label{BI heat}
\partial_t B+\nabla\times E=0,
\;\;\;E=\partial_1 H_\lambda(\nabla\times(\theta'(B)),B)
\end{split}
\end{equation}
by calling dissipative solution any admissible pair $(E,B)$
such that
\begin{equation}
\label{BI dissip}
\begin{split}
\frac{d}{dt}\int \theta(B)dx
+\int (L_\lambda(E,B)+H_\lambda(D,B))\le 0.
\\
{\rm{where\;\;}}D=\nabla\times (\theta'(B))
\\
H_\lambda(D,B)=
\sqrt{(\lambda^2+|B|^2)(1+|D|^2)-(D\cdot B)^2}.
\end{split}
\end{equation}
The Born-Infeld Lagrangian has two remarkable properties \cite{BI,Br}: as $\lambda\rightarrow \infty$,
we recover the classical Maxwell Lagrangian for electromagnetism
(as the limit of $(\lambda+L_\lambda)\lambda$); as  $\lambda\rightarrow 0$,
we get 
\begin{equation}
\label{BI cost}
H_0(D,B)=\sqrt{|B|^2(1+|D|^2)-(D\cdot B)^2},
\;\;\;L_0(E,B)=-\sqrt{|B|^2-|E|^2},\;\;E\cdot B=0,
\end{equation}
which, interestingly enough, includes the pointwise constraint $E\cdot B=0$.
This $exactly$
means there is a vector $v\in\ERRE^3$ such that $E=B\times v$
and $v$ can be defined, for instance, by setting 
$$
v=\frac{E\times B}{|B|^2}.
$$
So, the constraint $\partial_tB+\nabla\times E=0$ becomes
$$
\partial_t B+\nabla\times(B\times v)=0,
$$
or, equivalently,
$$
\partial_t B+\nabla\cdot(B\otimes v-v\otimes B)=0,
$$
which is exactly the "topology-preserving" equation (\ref{TE})
(with, here, a vector field $v$ that is a priori not divergence-free).
Let us now express $v$ in terms of $B$. We have
$$
E=\partial_1 H_0(D,B)=\frac{D|B|^2-(D\cdot B)B}{H_0(D,B)}
$$
where $D=\nabla\times(\theta'(B))$.
Thus
$$
v=\frac{E\times B}{|B|^2}=\frac{D\times B}{H_0(D,B)}=\frac{D\times B}{\sqrt{|B|^2+|D\times B|^2}}.
$$
So, in the case $\lambda=0$, equation (\ref{BI heat})
can be rephrased as
\begin{equation}
\begin{split}
\label{BI heat2}
\partial_t B+\nabla\cdot(B\otimes v-v\otimes B)=0,
\\
v=\frac{D\times B}{\sqrt{|B|^2+|D\times B|^2}}
\\
D=\nabla\times(\theta'(B)).
\end{split}
\end{equation}
The resulting system looks very much like the magnetic relaxation equation (\ref{MR})
(without divergence-free
constraint for $v$).
For instance, in the case  $\theta(B)=|B|$, we get 
$$
v=\nabla\cdot\frac{B\otimes B}{H_0},
$$
where $H_0=\sqrt{|B|^2+|D\times B|^2}$ and $D=\nabla\times (B/|B|)$.

%%%%%%%%%%%%%%%%%%%%%%%%%%%%%%%%%%%%%%%%%%%%%%
%%%%%%%%%%%%%%%%%%%%%%%%%%%%%%%%%%%%%%%%%%%%%%%%%%%%%%%%%%%%

\end{document}